\documentclass{elsarticle}
\usepackage{amssymb, color}
\usepackage{amsmath}
 \usepackage{pifont}
 \usepackage{caption}
 \usepackage[utf8]{inputenc}
\usepackage{setspace}
\makeatletter
\def\LaTeX{\leavevmode L\raise.42ex
    \hbox{\kern-.3em\size{\sf@size}{0pt}\selectfont A}\kern-.15em\TeX}
\makeatother

\newcommand{\BibTeX}{{\rm B\kern-.05em{\sc
          i\kern-.025emb}\kern-.08em\TeX}}

\makeatletter
\def\@currentlabel{2.1}\label{e:dispaa}
\def\@currentlabel{2.21}\label{e:dispau}
\def\@currentlabel{2.22}\label{e:dispav}
\def\@currentlabel{2.23}\label{e:dispaw}
\def\@currentlabel{2.24}\label{e:dispax}
\def\theequation{\thesection.\@arabic\c@equation}
\makeatother

\renewcommand{\theequation}{\arabic{section}.\arabic{equation}}

\newcommand{\R}{\mathbb R}

\newcommand{\e}{\epsilon}

\def \D{\Delta}

\newtheorem{thm}{Theorem} [section]
\newtheorem{lem}{Lemma} [section]
\newtheorem{prop}{Proposition} [section]

\newtheorem{rem}{Remark}[section]
\newtheorem{taggedtheoremx}{Theorem}
\newenvironment{taggedtheorem}[1]
 {\renewcommand\thetaggedtheoremx{#1}\taggedtheoremx}
 {\endtaggedtheoremx}

\renewcommand{\theequation}{\thesection.\arabic{equation}}
\renewcommand{\thesection}{\arabic{section}}
\renewcommand{\theequation}{\thesection.\arabic{equation}}
\let\ssection=\section\renewcommand{\section}{\setcounter{equation}{0}\ssection}

\begin{document}

\begin{frontmatter}
\title{Liouville type theorems for stable solutions of  elliptic system
involving the Grushin operator.}
\author[fd]{Foued Mtiri}
\ead{mtirifoued@yahoo.fr}
\address[fd]{ANLIG, UR13ES32, University of Tunis El-Manar, 2092 El Manar II, Tunisia.}
\begin{abstract}
We examine the degenerate elliptic system

$$-\Delta_{s} u = v^p, \quad -\Delta_{s} v= u^\theta,  \quad u,v>0 \quad\mbox{in }\; \mathbb{R}^N=\mathbb{R}^{N_1}\times \mathbb{R}^{N_2}, \quad\mbox{where }\;\;\;\; s \geq 0\;\; \mbox{and} \;\;p,\theta >0.$$
We prove that the system  has no smooth stable solution provided   $p,\theta >0$ and $N_s< 2 + \alpha + \beta,$ where
 $$\alpha = \frac{2(p+1)}{p\theta - 1} \quad\mbox{and} \quad \beta = \frac{2(\theta +1)}{p\theta - 1}.$$
 This result is an extension of some result in \cite{ MY}. In particular, we establish a new the integral estimate  for $u$ and $v$ \;(see Proposition 1.1), which is crucial to deal with the case $0 < p < 1.$
\end{abstract}
\end{frontmatter}
 \section{Introduction}
\setcounter{equation}{0}
We start by noting that throughout this article,  $N_s:=N_1+(1+s)N_2$ is called the
homogeneous dimension associated to the Grushin operator:
$$\Delta_{s}=\Delta_{x}+|x|^{2s}\Delta_{y},$$ where $s \geq 0,$ and
  $$\Delta_{x}:=\sum_{i=1}^{N_{1}}\frac{\partial^{2}}{\partial x_{i}^{2}},\quad\mbox{and} \;\; \Delta_{y}:=\sum_{j=1}^{N_{2}}\frac{\partial^{2}}{\partial y_{j}^{2}},$$ are Laplace operators with respect to $x\in \mathbb{R}^{N_1},$ $y\in \mathbb{R}^{N_2}$ and   $|x|^{2s}=\left(\sum \limits_{i=1}^{N_{1}}x_{i}^{2}\right)^{s}$.
  \medskip

In this paper, we are  interested in the classification of stable solutions to the following degenerate elliptic system
\begin{align}\label{1.1}
-\Delta_{s} u = v^p, \quad -\Delta_{s} v= u^\theta,  \quad u,v>0\quad\mbox{in }\; \mathbb{R}^N=\mathbb{R}^{N_1}\times \mathbb{R}^{N_2}, \quad\mbox{where }\;  p,\theta >0.
\end{align}
\medskip

 In the case  $s = 0,$ the system \eqref{1.1} becomes
\begin{align}\label{1.2}
-\Delta u = v^p, \quad-\Delta v= u^\theta,\quad u,v>0\quad\mbox{in }\; \mathbb{R}^N, \quad\mbox{where }\; p,\theta >0.
\end{align}
It is well studied and there
are many interesting results on classifying the solutions of this system for various  $p$ and $ \theta.$
\medskip

In the last decade, increased attention has focused on the system \eqref{1.1}.
 We note that the Sobolev critical hyperbola
\begin{align*}
\frac{1}{p+1} + \frac{1}{\theta +1} = \frac{N-2}{N} ,
\end{align*}
 which is introduced independently  by  Mitidieri  \cite{em} and
Van der Vorst \cite{van},  plays a crucial role in the analysis of \eqref{1.1}.  Among the main interests, the Lane–Emden conjecture says that there has no positive classical solution if and only if the pair $(p, \theta)$ lies below the Sobolev critical hyperbola, i.e.,

\medskip
 \textbf{ Conjecture.}  Let $p, \theta > 0$. If the pair $(p, \theta)$ is subcritical i.e., if  $(p, \theta)$ satisfies  \begin{align}\label{LE}
 \frac{1}{p+1} + \frac{1}{\theta +1} > \frac{N-2}{N},
\end{align}
  Then system \eqref{1.2} has no positive classical solutions.
\medskip

The full Lane-Emden conjecture is still open. Only partial results for nonexistence are known, and many researchers have made contribution in pushing the progress forward. We shall briefly present some
important recent developments of this conjecture
 \medskip

The case of radial solutions was solved by Mitidieri \cite{em1} and Serrin-Zou \cite{Zs} constructed positive radial solutions on and
above the critical hyperbola, i.e
 \begin{align*}
\frac{1}{p+1} + \frac{1}{\theta +1} \leq \frac{N-2}{N},
\end{align*}

  which means that the nonexistence theorem is optimal for radial solutions. For non-radial solutions of the Lane-Emden
system,  certain Liouville theorems are known. Denote the scaling exponents of system \eqref{1.2} by
 \begin{align}\label{ab}
  \alpha = \frac{2(p+1)}{p\theta - 1}, \quad \beta = \frac{2(\theta +1)}{p\theta - 1}, \quad p\theta > 1.
 \end{align}
 Then subcritical condition \eqref{LE} is equivalent to

 \begin{align}
 \label{LEbis} N < 2 + \alpha + \beta, \quad\mbox{for } \quad p\theta > 1.
\end{align}
 There are no positive classical super-solutions to \eqref{1.2} if
 \begin{align*}
p\theta \leq 1 \quad \mbox{or } \quad p\theta > 1 \quad \mbox{and } \quad \max\,\left\{\alpha, \beta\right\}>N-2,
\end{align*}
 see \cite{em1, ST, ZsS}. This implies the conjecture for $N = 1,2.$  Also, in \cite{bll} the conjecture is true for
  \begin{align*}
\min\,\left\{\alpha, \beta\right\}>\frac{N-2}{2} \quad \mbox{with } \quad (\alpha, \beta)\neq(\frac{N-2}{2},\frac{N-2}{2}).
\end{align*}
\smallskip

 In dimension $N = 3,$  Serrin-Zou \cite{ZsS} gave a proof for the nonexistence of polynomially bounded
solutions, an assumption that was relaxed later by Polácik, Quittner and Souplet \cite{pqs}.  More recently,  in dimension $N = 4 $, the conjecture was settled completely by Souplet \cite{phs} who provided, in dimensions $N\geq5$, a more restrictive new region for the exponents $(p, \theta)$ that insures that  if $(p, \theta)$ satisfies \eqref{LEbis}, along with $ \max\,\left\{\alpha, \beta\right\}>N-3,$ then system \eqref{1.2} has no positive classical solutions.

\medskip

On the other direction, the Liouville type theorem for the class of stable solutions
for both  system has attracted the attention of many mathematicians.  We refer \cite{cow, HU, Hfh, MY} for Lane–Emden system. The nonexistence
of positive stable classical solutions of \eqref{1.2} was examined in  \cite{cow}. In particular,  Cowan \cite{cow} proved that if $N\leq 10$, \eqref{1.2} has no stable solution for any $2\leq p\leq \theta$. This result was then  extended in \cite{ HU} to the Lane–Emden system with weights. Among other things, we also mention that Hajlaoui et.al. \cite{Hfh} improved the previous works \cite{cow, HU} and mainly obtained a new comparison property which is key to deal with the case $1 < p \leq \frac{4}{3}.$
\smallskip

 A new approach  which is independently obtained
in \cite{ MY}, allows the authors to prove  the following
Liouville theorem for classical stable  solutions of \eqref{1.2} for any $p, \theta > 0,$  satisfying \eqref{LEbis}.

\begin{taggedtheorem}{A}
If $p, \theta > 0$ satisfies \eqref{LE}, then \eqref{1.2} has no smooth stable solution.
\end{taggedtheorem}
In the case $0 < p < 1,$  the main tools  are based  on  the relationship between the stability for
the Lane-Emden system \eqref{1.2} and   the
stability  for the a fourth order problem, called the $m$-biharmonic equation  $$\Delta^{2}_{m} u := \Delta (|\Delta u|^{m-2}\Delta u) = u^\theta\quad\mbox{in } \quad \R^N .$$
  The case $p>1$ was handled by  the results of \cite{Hfh}. The main motivation was to provide a proof of the Lane-Emden conjecture dealing with positive stable solutions.
 \medskip

Coming back to the Lane-Emden system \eqref{1.1} for the general case  $s \geq 0,$  the Liouville property is less understood and is more complicated to deal with than $s = 0,$ because the the operator $\Delta_{s}$ no longer has symmetry and it degenerates on the manifold $\{0\}\times \mathbb{R}^{N_2}$   which causes further mathematical difficulties.
\smallskip

 In  \cite{DP} the author  extended some of Cowan results  \cite{cow},  in order to prove the nonexistence of positive stable  solutions of \eqref{1.1} with $\frac{4}{3}< p \leq \theta$  and $N_{s} \leq 10,$ or $1<p\leq \min(\frac{4}{3}, \theta)$ with additional assumption. The main idea used in \cite{DP} is a combination of stability inequality, comparison
principle and bootstrap argument. After that, this idea was exploited by Mtiri \cite{foi}, the author has
obtained a  some Liouville theorems for  stable solutions of \eqref{1.1}, see Theorem {\bf B} below. This results  improves the bound given in \cite{DP}.
\medskip

 Our main objective is to classify the stable solution of \eqref{1.1} for any $p, \theta > 0,$  and  the general case  $s \geq 0.$ In order to state our results more accurately,  we define  the notion of stability where we consider a general system given by
\begin{align}\label{1.222}
 -\Delta_{s} u = f(x,v),\quad -\Delta_{s} v= g(x,u)\;\; \mbox{in $K$, a bounded regular domain }\; \subset \R^N,
\end{align}
where $f,g\in C^1(K \times \mathbb{R}).$  Following Montenegro \cite{MO}, a smooth solution $(u,v)$ of \eqref{1.222} is said to
be stable in $K$ if the following eigenvalue problem
\begin{align*}
 -\Delta \xi = f_{v}(x,v)\zeta +\eta \xi, \quad -\Delta \zeta = g_{u}(x,u)\xi + \eta\zeta \quad \mbox{in }\, K
\end{align*}
has a nonnegative eigenvalue $\eta$, with a positive smooth eigenfunctions pair $(\xi, \zeta)$.

\medskip

The main result in this paper is the following

\begin{thm}
\label{main4}
If $p, \theta > 0$  and $N_{\alpha}$ satisfying $ N_{s} < 2 + \alpha + \beta$, then \eqref{1.1} has no smooth stable solution.
\end{thm}

\medskip
To prove Theorem \ref{main4}, we borrow crucially the idea in \cite{ MY} who established Theorem {\bf A}. Without loss of generality, we consider only $\theta > p > 0$ and $p \theta > 1$. As we will see soon, the $\theta > p \geq 1$ case can be handled by
the results in \cite{foi}, so our main concern is the case
$$\theta > \frac{1}{p} > 1 > p > 0.$$

\medskip

 The main difficulty arises from the fact that there has no works  in literature dealing with stable solutions (radial or not) for the $m$-biharmonic–type Grushin equation. Consequently, it is difficult to use the  technique  developed in \cite{ MY}. We overcome this difficulty, we shall derive new  the important   integral estimate  for $u$ and $v$ which is crucial to deal with the case $0 < p < 1.$  The following Proposition transforms our notion of a stable solution of \eqref{1.1} into an inequality which allows the use of arbitrary test functions.

\begin{prop}\label{p12bis}
\label{l.1} Let $(u,v)$ denote a stable solution of \eqref{1.1}  with $0 < p < 1.$ Then
\begin{align}
\label{stb}
 \theta p \int_{\mathbb{R}^N}  u^{\theta-1}\gamma^{2}dxdy \leq  \int_{\mathbb{R}^N}v^{1-p}|\D_{s}\gamma|^{2} dxdy,
\end{align}
for all $\gamma\in C_c^{\infty}(\mathbb{R}^N).$
\end{prop}

\begin{rem}
\begin{itemize}
\item If $s=0$, we obtain a similar result in \cite{ MY} see Theorem {\bf A}.
\item Our approach is more easier than those developed in \cite{ MY},  for $p, \theta > 0$. To the best of our knowledge, no general Liouville type result was known for stable solution of \eqref{1.1} with $0 < p < 1.$
    \item  We note also that the method used in the present paper can be applied to study the
weighted systems, and to more general class of degenerate operator, such as the $\Delta_{s}$ operator (see \cite{ nb, nb9}) of the form

$$\Delta_s:=\sum_{j=1}^N s^2_j\Delta_{x^{(j)}}\;\;\;  s:=(s_1, \dots, s_N): \mathbb{R}^N\to \mathbb{R^{N}},$$

where $s_i: \mathbb{R}^N\to \mathbb{R}, \quad i=1,\dots, N,$  are nonnegative continuous functions satisfying some properties such that $\Delta_s$ homogeneous of degree two with respect to a group dilation in $\mathbb{R}^N.$
\end{itemize}
\end{rem}

\medskip\noindent
 \section{Liouville Type Theorem.}

 This section is devoted to the proof of Theorem \ref{main4}. For convenience, we
always denote by $C$ a generic constant whose concrete values may change
from line to line or even in the same line. If this constant depends on an
arbitrary small number $\epsilon$, then we may denote it by $C_{\epsilon}$. We also use Young inequality in the form $ab\leq \epsilon a^{q}+C_{\epsilon} b^{q'}$ for $q, q'>1$ satisfying $\frac{1}{q}+\frac{1}{q'}=1$

\medskip
\setcounter{equation}{0}

As mentioned before, we need only to consider the case $\theta > p$ and $p \theta > 1$. We split the proof into two cases:
$\theta > p \geq 1$ and $\theta > p^{-1} > 1 > p > 0$.

\medskip
\subsection{The case $\theta> p^{-1} > 1>p > 0$. }
\medskip

Here we handle the case  $0<p<1.$ We begin with  proving the integral estimate  for $u$ and $v$

\medskip
\noindent{\bf  proof of Proposition \ref{p12bis}.} By the definition of stability, there exist smooth positive functions $\xi$, $\zeta$ and $\eta \geq 0$ such that
$$-\Delta_{s} \xi = pv^{p-1}\zeta + \eta \xi, \; \; -\Delta_{s}\zeta = \theta u^{\theta -1}\xi + \eta\zeta \quad \mbox{in }\; \mathbb{R}^N=\mathbb{R}^{N_1}\times \mathbb{R}^{N_2}.$$
Using $(\xi, \zeta)$ as super-solution, $(\min_{\overline\Omega}\xi, \min_{\overline\Omega}\zeta)$ as sub-solution, and the standard monotone iterations, we can claim that
there exist positive smooth functions $\varphi$, $\chi$ verifying
\begin{align*}
 -\Delta_{s} \varphi = p v^{p-1}\chi, \quad -\Delta_{s} \chi = \theta u^{\theta-1}\varphi\quad \mbox{in }\, \mathbb{R}^N=\mathbb{R}^{N_1}\times \mathbb{R}^{N_2}.
\end{align*}
Therefore, we have
\begin{align*}
\theta u^{\theta-1}\varphi=\Delta_{s}\left(\frac{1}{p} v^{1-p} \Delta_{s} \varphi\right) \quad \mbox{in}\;\;\mathbb{R}^N=\mathbb{R}^{N_1}\times \mathbb{R}^{N_2}.
\end{align*}
Let $\gamma \in C_c^2(\mathbb{R}^N)$. Multiplying the above equation by $\gamma^{2}\varphi^{-1}$ and
integrating by parts, there holds
\begin{align}\label{0.255}
\begin{split}
 &\int_{\mathbb{R}^N} \theta u^{\theta-1}\gamma^{2}dxdy\\
  & = \frac{1}{p}\int_{\mathbb{R}^N} v^{1-p} \Delta_{s} \varphi\Delta_{s}(\gamma^{2}\varphi^{-1})dxdy\\
&= \frac{1}{p} \int_{\mathbb{R}^N}v^{1-p} \Delta_{s} \varphi\left(-4\gamma\frac{\nabla_{s} \varphi\cdot\nabla_{s}\gamma}{\varphi^{2}}+\frac{2|\nabla_{s}\gamma |^2}{\varphi}+\frac{2\gamma\D_{s}\gamma}{\varphi}+\frac{2\gamma^{2}|\nabla_{s}\varphi |^2}{\varphi^{3}} -\frac{\gamma^{2}\D_{s}\varphi}{\varphi^{2}}\right) dxdy.
\end{split}
\end{align}
Using Cauchy-Schwarz's inequality and the fact that $-\Delta_{s} \varphi >0,$ we get
\begin{align}\label{2.La}
\begin{split}
&\left| -4\int_{\mathbb{R}^N}\frac{v^{1-p}}{p} \Delta_{s} \varphi\frac{\nabla_{s} \varphi\cdot\nabla_{s}\gamma}{\varphi^{2}}\gamma dxdy\right|\\
 &\leq -2\int_{\mathbb{R}^N}\frac{v^{1-p}}{p} \Delta_{s} \varphi\frac{|\nabla_{s}\gamma |^2}{\varphi} dxdy -2\int_{\mathbb{R}^N}\frac{v^{1-p}}{p} \Delta_{s} \varphi\frac{\gamma^{2}|\nabla_{s}\varphi |^2}{\varphi^{3}} dxy.
\end{split}
\end{align}
Combining \eqref{0.255} and  \eqref{2.La}, one obtains, using again the Cauchy-Schwartz inequality,
\begin{align*}
 \int_{\mathbb{R}^N} \theta u^{\theta-1}\gamma^{2}dxdy &
\leq \frac{2}{p} \int_{\mathbb{R}^N}v^{1-p} \Delta_{s} \varphi\frac{\gamma\D_{s}\gamma}{\varphi} dxdy-\frac{1}{p} \int_{\mathbb{R}^N}v^{1-p}\frac{(\Delta_{s} \varphi)^{2}}{\varphi^{2}}\gamma^{2} dxdy\\
& \leq  \frac{1}{p}\int_{\mathbb{R}^N}v^{1-p}\frac{(\Delta_{s} \varphi)^{2}}{\varphi^{2}}\gamma^{2} dxdy + \frac{1}{p}\int_{\mathbb{R}^N}v^{1-p}(\D_{s}\gamma)^{2} dxdy -\frac{1}{p} \int_{\mathbb{R}^N}v^{1-p}\frac{(\Delta_{s} \varphi)^{2}}{\varphi^{2}}\gamma^{2} dxdy\\
& = \frac{1}{p} \int_{\mathbb{R}^N}v^{1-p}(\D_{s}\gamma)^{2} dxdy.
\end{align*}
 The proof is completed.\qed

\medskip

 As a consequence of Proposition \ref{p12bis}, we derive immediately the following integral estimate for  stable solution $(u,v)$ of system \eqref{1.1}, in the case $p \in (0, 1)$ and $p\theta > 1,$  which is a crucial tool in our approach.
  \begin{lem}
\label{lemnewBN} Let $(u,v)$  be  a stable solution of \eqref{1.1}, with  $p \in (0, 1)$. Then, for any integer  $$k\geq \max \left(\frac{1+p}{2p}, \frac{\beta(p+1)}{4}, \frac{\alpha(\theta+1)}{4}\right),$$
there exists a positive constant  $ C=C(N, \epsilon, p, k)$ such that for any $\zeta\in C_c^2(\mathbb{R}^N)$ satisfying $0 \leq\zeta \leq 1$,
\begin{align}\label{new7}
\begin{split}
\int_{\mathbb{R}^N}v^{p+1}\zeta^{4k}dxdy +\int_{\mathbb{R}^N} u^{\theta+1}\zeta^{4k} dxdy
\leq C\left[\int_{\mathbb{R}^N}\Big(|\D_{s} \zeta|^{\frac{p+1}{p}}+|\nabla_{s} \zeta|^{\frac{2(p+1)}{p}} +|\nabla_{s}^{2} \zeta|^{\frac{p+1}{p}}\Big)^{\frac{p\beta}{2}} dxdy\right].
 \end{split}
\end{align}
Here $\alpha$ and $\beta$ are defined in \eqref{ab}.
\end{lem}
\noindent{\bf Proof.}  First, for any  $\epsilon \in (0, 1)$ and $\eta \in C^2(\R^N)$, there holds
\begin{align}
\label{newesTt47}
  \begin{split}
\int_{\mathbb{R}^N}v^{1-p} [\Delta_{s} (u\eta)]^2 dxdy & = \int_{\mathbb{R}^N}v^{1-p}\left(u \D_{s} \eta+2\nabla_{s} u\nabla_{s} \eta + \eta\D_{s} u\right)^{2} dxdy\\
& \leq \left(1+C\epsilon\right)\int_{\mathbb{R}^N}v^{p+1} {\eta}^2 dxdy + \frac{C}{\e}\int_{\mathbb{R}^N}v^{1-p} u^2 |\D_{s}\eta|^2 dxdy \\
  &+\frac{C}{\epsilon} \int_{\mathbb{R}^N} v^{1-p}|\nabla_{s} u|^2|\nabla_{s} \eta|^2 dxdy.
    \end{split}
   \end{align}
    Take $\eta = \zeta^{2k}$ with $k \geq \frac{1}{p}+1 > 2$ and  $0 \leq\zeta \leq 1.$  Apply Young's inequality, we get
     \begin{align*}
  \int_{\mathbb{R}^N} v^{1-p} u^2|\D_{s} (\zeta^{2k})|^2 dxdy & \leq C_{k}\int_{\mathbb{R}^N}v^{1-p} u^2\left(|\D_{s} \zeta|^{2}+|\nabla_{s} \zeta|^4\right)\zeta^{\frac{4k(1-p)}{1+p}+\frac{2(4pk-2(1+p)}{1+p}} dxdy\\
& \leq C_{p, k}\epsilon^{2} \int_{\mathbb{R}^N}v^{p+1} \zeta^{4k} dxdy +
  C_{\epsilon, k, p}\int_{\mathbb{R}^N} u^{\frac{1+p}{p}} \left(|\D_{s} \zeta|^{2}+|\nabla_{s} \zeta|^4\right)^{\frac{p+1}{2p}}\zeta^{4k-2(\frac{1}{p}+1)} dxdy,
\end{align*}
and
  \begin{align*}
  \int_{\mathbb{R}^N} v^{1-p}|\nabla_{s} u|^2|\nabla_{s} (\zeta^{2k})|^2 dxdy & =4 k^{2}\int_{\mathbb{R}^N} v^{1-p}|\nabla_{s} u|^2|\nabla_{s} \zeta|^2\zeta^{\frac{4k(1-p)}{1+p}+\frac{2(4pk-(1+p)}{1+p}} dxdydy\\
& \leq   C_{p, k}\epsilon^{2} \int_{\mathbb{R}^N}v^{p+1} \zeta^{4k} dxdy +
  \frac{C_{p,k}}{\epsilon^{2}}\int_{\mathbb{R}^N} |\nabla_{s} u|^{\frac{p+1}{p}} |\nabla_{s} \zeta|^{\frac{p+1}{p}}\zeta^{4k-\frac{p+1}{p}} dxdy.
\end{align*}
 Combining all these inequalities with \eqref{newesTt47},  we get the following estimates

 \begin{align}\label{newest4}
  \begin{split}
&\int_{\mathbb{R}^N}v^{1-p} [\Delta_{s} (u\zeta^{2k})]^2 dxdy\\
  &\leq \left(1+C_{p, k}\epsilon\right) \int_{\mathbb{R}^N}v^{p+1} \zeta^{4k} dxdy+\frac{C_{p,k}}{\epsilon^{3}}\int_{\mathbb{R}^N} |\nabla_{s} u|^{\frac{p+1}{p}} |\nabla_{s} \zeta|^{\frac{p+1}{p}}\zeta^{4k-\frac{p+1}{p}} dxdy\\
  &+ C_{\epsilon, k, p}\int_{\mathbb{R}^N} u^{\frac{1+p}{p}} \left(|\D_{s} \zeta|^{2}+|\nabla_{s} \zeta|^4\right)^{\frac{p+1}{2p}}\zeta^{4k-2(\frac{1}{p}+1)} dxdy.
  \end{split}
   \end{align}

\medskip

We will use also the following lemma
 \begin{lem}\label{l.2.7a}
  Let $k \geq \frac{1+p}{2p} > 1,$ there exists a positive constant $C=C_{N, \e, p, k}>0$ such that for any $\epsilon >0,$ and $\zeta \in C_c^{\infty}(\mathbb{R}^N)$ with $0 \leq\zeta \leq 1$, there holds

 \begin{align*}
&\;\int_{\mathbb{R}^N} |\nabla_{s} u|^{\frac{p+1}{p}} |\nabla_{s} \zeta|^{\frac{p+1}{p}}\zeta^{4k-\frac{p+1}{p}} dxdy \\
 &\;\leq \epsilon^{4}\int_{\mathbb{R}^N}v^{p+1}\zeta^{4k}dxdy +C \int_{\mathbb{R}^N}u^{\frac{p+1}{p}}\left(|\nabla_{s} \zeta|^{\frac{2(p+1)}{p}} +|\nabla_{s}^{2} \zeta|^{\frac{p+1}{p}}\right)\zeta^{4k-\frac{2(p+1)}{p}} dxdy.
\end{align*}
\end{lem}
 \noindent{\bf Proof.}
Since the simple calculation implies that
\begin{align}
  \label{new1}
  \begin{split}
 \int_{\mathbb{R}^N} |\nabla_{s} u|^{\frac{p+1}{p}} |\nabla_{s} \zeta|^{\frac{p+1}{p}}\zeta^{4k-\frac{p+1}{p}} dxdy = &\;\int_{\mathbb{R}^N} \nabla_{s} u\cdot\nabla_{s} u|\nabla_{s} u|^{\frac{p+1}{p}-2} |\nabla_{s} \zeta|^{\frac{p+1}{p}}\zeta^{4k-\frac{p+1}{p}} dxdy\\
 =& \; - \int_{\mathbb{R}^N}  \mathrm{div}_{s} \left(\nabla_{s} u|\nabla_{s} u|^{\frac{p+1}{p}-2}\right) u|\nabla_{s} \zeta|^{\frac{p+1}{p}}\zeta^{4k-\frac{p+1}{p}} dxdy\\
 &\; -\int_{\mathbb{R}^N} u|\nabla_{s} u|^{\frac{p+1}{p}-2}\nabla_{s} u \cdot\nabla_{s}\left(|\nabla_{s} \zeta|^{\frac{p+1}{p}}\zeta^{4k-\frac{p+1}{p}}\right)dxdy\\
 := & \; Q_1 + Q_2,
 \end{split}
\end{align}
 where $ \mathrm{div}_{s}= \mathrm{div}_{x}+|x|^{s}\mathrm{div}_{y}.$ hence the first term on the right hand side of \eqref{new1}, can be estimated as
 \begin{align*}
  Q_1 & = -\left(\frac{p+1}{p}-2\right) \int_{\mathbb{R}^N}  u |\nabla_{s} u|^{\frac{p+1}{p}-4}|\nabla_{s} \zeta|^{\frac{p+1}{p}}\nabla_{s}^{2}u (\nabla_{s} u, \nabla_{s} u)\zeta^{4k-\frac{p+1}{p}} dxdy \\
 &\;+\int_{\mathbb{R}^N}u v^p |\nabla_{s} u|^{\frac{p+1}{p}-2}|\nabla_{s} \zeta|^{\frac{p+1}{p}}\zeta^{4k-\frac{p+1}{p}} dxdy &\;\\
 & \leq C_p\int_{\mathbb{R}^N}  u| \nabla_{s}^{2} u||\nabla_{s} u|^{\frac{p+1}{p}-2}|\nabla_{s} \zeta|^{\frac{p+1}{p}}\zeta^{4k-\frac{p+1}{p}} dxdy  +\int_{\mathbb{R}^N}u v^p |\nabla_{s} u|^{\frac{p+1}{p}-2}|\nabla_{s} \zeta|^{\frac{p+1}{p}}\zeta^{4k-\frac{p+1}{p}} dxdy .
\end{align*}
 Observe that $$ 4k-\frac{p+1}{p}=\left(\frac{4k-\frac{p+1}{p}}{\frac{p+1}{p}}\right)\left(\frac{p+1}{p}-2\right)+2\left(\frac{4k-\frac{p+1}{p}}{\frac{p+1}{p}}\right),$$

  and  applying Young's inequality, there holds, for any $\epsilon > 0$,

\begin{align}\label{0.25512vc}
 \begin{split}
 &\int_{\mathbb{R}^N}  u v^p |\nabla_{s} u|^{\frac{p+1}{p}-2}|\nabla_{s} \zeta|^{\frac{p+1}{p}-2+2}\zeta^{4k-\frac{p+1}{p}} dxdy\\
 \leq & \;C_{ \epsilon, p}\int_{\mathbb{R}^N} \left(u v^p\right)^{\frac{p+1}{2p}}|\nabla_{s} \zeta|^{\frac{p+1}{p}}\zeta^{4k-\frac{p+1}{p}} dxdy + C_{p}\epsilon\int_{\mathbb{R}^N}|\nabla_{s} u|^{\frac{p+1}{p}}|\nabla_{s} \zeta|^{\frac{p+1}{p}}\zeta^{4k-\frac{p+1}{p}} dxdy\\
 \leq &\;C_{\epsilon, p}\int_{\mathbb{R}^N}  u^{\frac{p+1}{p}}|\nabla_{s} \zeta|^{\frac{2(p+1)}{p}}\zeta^{4k-\frac{2(p+1)}{p}} dxdy+ C_{p}\epsilon^{4}\int_{\mathbb{R}^N}v^{p+1}\zeta^{4k}dxdy\\
 &+ C_{p}\epsilon\int_{\mathbb{R}^N}|\nabla_{s} u|^{\frac{p+1}{p}}|\nabla_{s} \zeta|^{\frac{p+1}{p}}\zeta^{4k-\frac{p+1}{p}} dxdy.
\end{split}
\end{align}

and
 \begin{align}\label{0.25512}
 \begin{split}
 & \;C_p\int_{\mathbb{R}^N}  u| \nabla_{s}^{2} u||\nabla_{s} u|^{\frac{p+1}{p}-2}|\nabla_{s} \zeta|^{\frac{p+1}{p}-2+2}\zeta^{4k-\frac{p+1}{p}} dxdy\\
 \leq &\;  C_{\epsilon, p}\int_{\mathbb{R}^N}\left(u | \nabla_{s}^{2} u|\right)^{\frac{p+1}{2p}}|\nabla_{s} \zeta|^{\frac{p+1}{p}}\zeta^{4k-\frac{p+1}{p}} dxdy  + + C_{p}\epsilon\int_{\mathbb{R}^N}|\nabla_{s} u|^{\frac{p+1}{p}}|\nabla_{s} \zeta|^{\frac{p+1}{p}}\zeta^{4k-\frac{p+1}{p}} dxdy\\
 \leq &\;C_{\epsilon, p}\int_{\mathbb{R}^N}  u^{\frac{p+1}{p}}|\nabla_{s} \zeta|^{\frac{2(p+1)}{p}}\zeta^{4k-\frac{2(p+1)}{p}} dxdy+ C_{p}\epsilon^{4}\int_{\mathbb{R}^N}| \nabla_{s}^{2} u|^{\frac{p+1}{p}}\zeta^{4k}dxdy\\
 &+ C_{p}\epsilon\int_{\mathbb{R}^N}|\nabla_{s} u|^{\frac{p+1}{p}}|\nabla_{s} \zeta|^{\frac{p+1}{p}}\zeta^{4k-\frac{p+1}{p}} dxdy.
\end{split}
\end{align}
Now we shall estimate the integral
$$\int_{\mathbb{R}^N}| \nabla_{s}^{2} u|^{\frac{p+1}{p}}\zeta^{4k}dxdy.$$

For the  last inequality, we have used the following estimates from [\cite{ nb99}, Proof of  Theorem 2.6]

  $$\|\nabla_{x_{i}x_{j}}u\|_{L^{2}(\R^{N_{1}})} \leq C(N_{1})\|\D_{x}u\|_{L^{2}(\R^{N_{1}})}, \quad \mbox{for}\quad 1\leq i,j\leq N_{1},$$

and
$$\|\nabla_{y_{h}y_{k}}u\|_{L^{2}(\R^{N_{2}})} \leq C(N_{2})\|\D_{y}u\|_{L^{2}(\R^{N_{2}})}, \quad \mbox{for}\quad 1\leq h,k\leq N_{2}.$$

Let $\psi^{r}\in C_0^{\infty}(\mathbb{R}^N),$ with  $r>2.$ By direct calculations, we get

  \begin{align*}
|\nabla_{s}^2(u)|\psi^{r}
 \leq & \; C_{r}\Big[u\left(|\nabla_{s} \psi|^{2}\psi^{r-2}+|\nabla_{s}^{2} \psi|\psi^{r-1}\right)+ |\nabla_{s} u||\nabla_{s} \psi|\psi^{r-1}+ |\nabla_{s}^2(u\psi^{r})|\Big].
\end{align*}

Consider $\psi = \zeta,$ and  $r = \frac{4kp}{p+1}\geq 2, $ so that $k\geq\frac{p+1}{2p}.$ We can claim, for any $0 \leq\zeta \leq 1,$   there exists $C_{p,k} > 0$ such that
 \begin{align}\label{0.2551}
 \begin{split}
 & \;\int_{\mathbb{R}^N}| \nabla_{s}^{2} u|^{\frac{p+1}{p}}\zeta^{4k}dxdy\\
 & \;\leq C_{p,k}\int_{\mathbb{R}^N}|\nabla_{s}^2(u\zeta^{\frac{4kp}{p+1}})|^{\frac{p+1}{p}} dxdy +C_{p,k}\int_{\mathbb{R}^N}|\nabla_{s} u|^{\frac{p+1}{p}}|\nabla_{s} \zeta|^{\frac{p+1}{p}}\zeta^{4k-\frac{p+1}{p}} dxdy\\
 &\;+ C_{p,k}\int_{\mathbb{R}^N}u^{\frac{p+1}{p}}\left(|\nabla_{s} \zeta|^{\frac{2(p+1)}{p}} +|\nabla_{s}^{2} \zeta|^{\frac{p+1}{p}}\right)\zeta^{4k-\frac{2(p+1)}{p}} dxdy.
 \end{split}
\end{align}

Clearly $u\zeta \in H_{0}^{2}(\R^N),$ so from a standard approximation and scaling argument imply then
\begin{align}\label{0.j255}
 \begin{split}
\int_{\mathbb{R}^N}|\nabla_{s}^2(u\zeta^{\frac{4kp}{p+1}})|^{\frac{p+1}{p}} dxdy & \leq\;
C_{N, p} \int_{\mathbb{R}^N}|\D_{s}(u\zeta^{\frac{4kp}{p+1}})|^{\frac{p+1}{p}} dxdy\\
& \leq\; C_{N, p,k}\int_{\mathbb{R}^N}|\nabla_{s} u|^{\frac{p+1}{p}}|\nabla_{s} \zeta|^{\frac{p+1}{p}}\zeta^{4k-\frac{p+1}{p}} dxdy+ C_{N, p}\int_{\mathbb{R}^N}v^{p+1}\zeta^{4k}dxdy\\
  &+\;C_{N, p,k}\int_{\mathbb{R}^N}u^{\frac{p+1}{p}}\left(|\nabla_{s} \zeta|^{\frac{2(p+1)}{p}} +|\nabla_{s}^{2} \zeta|^{\frac{p+1}{p}}\right)\zeta^{4k-\frac{2(p+1)}{p}} dxdy.
\end{split}
\end{align}

Combining \eqref{0.25512vc}--\eqref{0.j255}, we obtain the estimate for the first left term in \eqref{new1}:

 \begin{align*}
Q_1 \leq &\;C_{N, p,k}\epsilon^{4}\int_{\mathbb{R}^N}v^{p+1}\zeta^{4k}dxdy +C_{N, p,k}\epsilon\int_{\mathbb{R}^N}|\nabla_{s} u|^{\frac{p+1}{p}}|\nabla_{s} \zeta|^{\frac{p+1}{p}}\zeta^{4k-\frac{p+1}{p}} dxdy\\
  &+\; C_{N, p,k}\int_{\mathbb{R}^N}u^{\frac{p+1}{p}}\left(|\nabla_{s} \zeta|^{\frac{2(p+1)}{p}} +|\nabla_{s}^{2} \zeta|^{\frac{p+1}{p}}\right)\zeta^{4k-\frac{2(p+1)}{p}} dxdy.
\end{align*}

Furthermore, by Young’s inequality,
 \begin{align*}
Q_2 = & - \frac{p+1}{p}\int_{\mathbb{R}^N} u\Big(|\nabla_{s} u||\nabla_{s} \zeta|\Big)^{\frac{p+1}{p}-2}\nabla_{s}^2\zeta(\nabla_{s}\zeta, \nabla_{s} u)\zeta^{4k-\frac{p+1}{p}} dxdy\\
& -(4k-\frac{p+1}{p})\int_{\mathbb{R}^N} u|\nabla_{s} u|^{\frac{p+1}{p}-2}|\nabla_{s} \zeta|^{\frac{p+1}{p}}(\nabla_{s} u\cdot \nabla_{s}\zeta)\zeta^{4k-\frac{p+1}{p}-1}dxdy\\
\leq &\; C_{ p,k} \int_{\mathbb{R}^N} u\Big(|\nabla_{s} u||\nabla_{s} \zeta|\Big)^{\frac{p+1}{p}-1}\left(|\nabla_{s} \zeta|^{2} +|\nabla_{s}^{2} \zeta|\right)\zeta^{4k-\frac{p+1}{p}-1}dxdy\\
\leq &\; C_{\epsilon, p,k} \int_{\mathbb{R}^N}u^{\frac{p+1}{p}}\left(|\nabla_{s} \zeta|^{\frac{2(p+1)}{p}} +|\nabla_{s}^{2} \zeta|^{\frac{p+1}{p}}\right)\zeta^{4k-\frac{2(p+1)}{p}} dxdy
+ \epsilon\int_{\mathbb{R}^N}|\nabla_{s} u|^{\frac{p+1}{p}}|\nabla_{s} \zeta|^{\frac{p+1}{p}}\zeta^{4k-\frac{p+1}{p}} dxdy.
\end{align*}

Combining the last tow inequality with \eqref{new1}, we get readily
 \begin{align*}
&\;(1-C_{N, p,k}\epsilon)\int_{\mathbb{R}^N}|\nabla_{s} u|^{\frac{p+1}{p}}|\nabla_{s} \zeta|^{\frac{p+1}{p}}\zeta^{4k-\frac{p+1}{p}} dxdy\\
& \leq C_{\epsilon, p, N ,k} \int_{\mathbb{R}^N}u^{\frac{p+1}{p}}\left(|\nabla_{s} \zeta|^{\frac{2(p+1)}{p}} +|\nabla_{s}^{2} \zeta|^{\frac{p+1}{p}}\right)\zeta^{4k-\frac{2(p+1)}{p}} dxdy\\
  &\;+C_{N, p,k}\epsilon^{4}\int_{\mathbb{R}^N}v^{p+1}\zeta^{4k}dxdy.
\end{align*}
Take $\epsilon$ small enough, the lemma  follows.\qed

\medskip

Now, using Lemma \ref{l.2.7a} and \eqref{newest4} , we obtain also
\begin{align}\label{nt4}
  \begin{split}
\int_{\mathbb{R}^N}v^{1-p} [\Delta_{s} (u\zeta^{2k})]^2 dxdy
  \leq& \; C\int_{\mathbb{R}^N}u^{\frac{p+1}{p}}\left(|\D_{s} \zeta|^{\frac{p+1}{p}}+|\nabla_{s} \zeta|^{\frac{2(p+1)}{p}} +|\nabla_{s}^{2} \zeta|^{\frac{p+1}{p}}\right)\zeta^{4k-\frac{2(p+1)}{p}} dxdy\\
  & \;+\left(1+C\epsilon\right) \int_{\mathbb{R}^N}v^{p+1}\zeta^{4k}dxdy.
  \end{split}
   \end{align}
Thanks to the approximation argument, the stability property \eqref{stb}, with $\gamma=u \zeta^{2k}$ where $\zeta \in C_0^{\infty}(\mathbb{R}^N),$ we get \begin{align}\label{f1}
 \begin{split}
& \;\theta p\int_{\mathbb{R}^N} u^{\theta+1}\zeta^{4k} dxdy -\left(1+ C\epsilon\right) \int_{\mathbb{R}^N}v^{p+1}\zeta^{4k}dxdy\\
\leq  & \;C\int_{\mathbb{R}^N}u^{\frac{p+1}{p}}\left(|\D_{s} \zeta|^{\frac{p+1}{p}}+|\nabla_{s} \zeta|^{\frac{2(p+1)}{p}} +|\nabla_{s}^{2} \zeta|^{\frac{p+1}{p}}\right)\zeta^{4k-\frac{2(p+1)}{p}} dxdy.
  \end{split}
\end{align}

Multiplying the equation $ -\Delta_{s} v= u^\theta$ by $u \zeta^{4k}$ and integrating by parts, there holds

\begin{align*}
 \int_{\mathbb{R}^N}v^{p+1}\zeta^{4k}dxdy-\int_{\mathbb{R}^N} u^{\theta+1}\zeta^{4k} dxdy \leq\int_{\mathbb{R}^N}u v
  |\D(\zeta^{4k})| dxdy
  + \int_{\mathbb{R}^N}v |\nabla_{s} u||\nabla_{s} (\zeta^{4k})| dxdy.
\end{align*}

Using Young's inequality and applying again Lemma \ref{l.2.7a}, we can conclude that for any $\epsilon>0$, there exists $C=C_{N, \e, p, k}>0$
 such that
\begin{align}\label{f12}
\begin{split}
& \;\left(1- C\epsilon\right) \int_{\mathbb{R}^N}\int_{\mathbb{R}^N}v^{p+1}\zeta^{4k}dxdy-\int_{\mathbb{R}^N} u^{\theta+1}\zeta^{4k} dxdy\\
\leq& \;C\int_{\mathbb{R}^N}u^{\frac{p+1}{p}}\left(|\D_{s} \zeta|^{\frac{p+1}{p}}+|\nabla_{s} \zeta|^{\frac{2(p+1)}{p}} +|\nabla_{s}^{2} \zeta|^{\frac{p+1}{p}}\right)\zeta^{4k-\frac{2(p+1)}{p}} dxdy.
\end{split}
\end{align}

Now,  multiplying \eqref{f12} by $\frac{1+ 2C\epsilon}{1- C\epsilon}$, adding it with \eqref{f1}, we get
\begin{align*}
& C\epsilon \int_{\mathbb{R}^N}v^{p+1}\zeta^{4k}dxdy + \left(p\theta-\frac{1+ 2C\epsilon}{1- C\epsilon}\right)\int_{\mathbb{R}^N} u^{\theta+1}\zeta^{4k} dxdy\\
\leq  & \; C\int_{\mathbb{R}^N}u^{\frac{p+1}{p}}\left(|\D_{s} \zeta|^{\frac{p+1}{p}}+|\nabla_{s} \zeta|^{\frac{2(p+1)}{p}} +|\nabla_{s}^{2} \zeta|^{\frac{p+1}{p}}\right)\zeta^{4k-\frac{2(p+1)}{p}} dxdy.
\end{align*}

 As $\theta > p^{-1}> 1.$ Fix  $0 <\epsilon < \frac{p\theta-1}{C( p\theta+ 2)}$, there holds
\begin{align}\label{f1xx2xx}
\begin{split}
& \int_{\mathbb{R}^N}v^{p+1}\zeta^{4k}dxdy +\int_{\mathbb{R}^N} u^{\theta+1}\zeta^{4k} dxdy\\
\leq  & \; C\int_{\mathbb{R}^N}u^{\frac{p+1}{p}}\left(|\D_{s} \zeta|^{\frac{p+1}{p}}+|\nabla_{s} \zeta|^{\frac{2(p+1)}{p}} +|\nabla_{s}^{2} \zeta|^{\frac{p+1}{p}}\right)\zeta^{4k-\frac{2(p+1)}{p}} dxdy.
\end{split}
\end{align}
Apply Young's inequality, we deduce then
\begin{align*}
& \; \int_{\mathbb{R}^N}v^{p+1}\zeta^{4k}dxdy +\int_{\mathbb{R}^N} u^{\theta+1}\zeta^{4k} dxdy\\
\leq  & \; C_{\epsilon}\int_{\mathbb{R}^N}\Big(|\D_{s} \zeta|^{\frac{p+1}{p}}+|\nabla_{s} \zeta|^{\frac{2(p+1)}{p}} +|\nabla_{s}^{2} \zeta|^{\frac{p+1}{p}}\Big)^{\frac{p\beta}{2}} dxdy+\epsilon' \int_{\mathbb{R}^N} |u|^{\theta+1}\zeta^{\frac{4kp(\theta+1)}{p+1}-2(\theta+1)}dxdy \\
\leq  & \; C_{\epsilon}\int_{\mathbb{R}^N}\Big(|\D_{s} \zeta|^{\frac{p+1}{p}}+|\nabla_{s} \zeta|^{\frac{2(p+1)}{p}} +|\nabla_{s}^{2} \zeta|^{\frac{p+1}{p}}\Big)^{\frac{p\beta}{2}} dxdy+\epsilon' \int_{\mathbb{R}^N} |u|^{\theta+1}\zeta^{4k}dxdy .
 \end{align*}

 Choose  $k\geq \frac{\beta(p+1)}{4}=\frac{\alpha(\theta+1)}{4} $  so that $4k\leq \frac{4kp(\theta+1)}{p+1}-2(\theta+1),$ for the last line. Take $\epsilon'$ small enough, the estimate \eqref{new7} is proved. \qed

\medskip

We are now in position to conclude. Choose $\chi $ a cut-off function in $ C_c^\infty\left(\mathbb{R}^N=\mathbb{R}^{N_1}\times\mathbb{R}^{N_2}, [o,1]\right),$ such that $$\chi=1 \quad \mbox{on} \quad B_{1}\times B_{2},\quad \mbox{and} \quad  \chi=0 \quad \mbox{outside } \quad B_{2}\times B_{2^{1+s}}.$$

For $R>0,$ put $\eta_{R}(x,y)=\chi(\frac{x}{R},\frac{y}{R^{1+s}}),$ it is easy to verify  that there exists $C >0$ independent of $R$ such that

$$|\nabla_{x}\eta_{R}|\leq \frac{C}{R}\quad \mbox{and} \quad  |\nabla_{y}\eta_{R}|\leq \frac{C}{R^{1+s}},$$

\smallskip

$$|\nabla^{2}_{x}\eta_{R}|+|\D_{x}\eta_{R}|\leq \frac{C}{R^{2}}\quad \mbox{and} \quad |\nabla^{2}_{y}\eta_{R}|+ |\D_{y}\eta_{R}|\leq \frac{C}{R^{2(1+s)}}.$$

Applying \eqref{new7} with $\zeta = \eta_{R}(x,y),$ we get
 \begin{align*}
\int_{ B_{R}\times B_{R^{1+s}}} u^{\theta+1}dxdy\leq \int_{\R^N} u^{\theta+1}\eta_{R}^{4k} dxdy \leq C R^{\frac{2N_{s}}{p\beta} -\frac{2(p+1)}{p}}.
\end{align*}

Since under our assumptions $ N_{s} < \frac{2(p+1)(\theta+1)}{p\theta - 1}:=2 + \alpha + \beta,$ the desired claim follows by letting $R \to \infty.$ \qed
\medskip

\subsection{ The case $\theta>p \geq 1$.}
\medskip

Let us recall a consequence of Theorem 1.1 (with $\alpha = 0$ there) in \cite{foi}.

\begin{taggedtheorem}{B}
\label{main3}
Let $x_0$ be the largest root of the polynomial
\begin{align}\label{newH}
H(x)=x^4 -p\theta\alpha\beta \left[4x^{2}-2(\alpha+\beta)x+1\right].
\end{align}
  \begin{enumerate}
\item
 If $\frac{4}{3}< p \leq \theta$ then \eqref{1.1} has no stable  solution if $N_{s}<2+2x_0.$

\item If $1<p\leq \min(\frac{4}{3}, \theta)$, then \eqref{1.1} has no bounded  stable solution, if
$$N_{s} < 2 + 2x_0\left[\frac{p}{2}+\frac{(2-p)(p \theta -1)}{(\theta+p-2)(\theta+1)}\right].$$
\end{enumerate}
\end{taggedtheorem}
\medskip

For $\theta > p \geq 1,$ we can proceed
similarly as the proof of Proposition in \cite{MY}, we sketch a proof here for the reader’s convenience.
Performing the change of variables $x=\frac{\beta}{2}z$ in \eqref{newH}, a direct computation shows that $H(x)=\left(\frac{\beta}{2}\right)^4L(z)$ where
$$L(z):=z^4-\frac{16p\theta(p+1)}{\theta+1}z^2+\frac{16p\theta(p+1)(p+\theta+2)}{(\theta+1)^2}z-\frac{16p\theta(p+1)^2}{(\theta+1)^2}.$$
Denote by $z_0$  the largest root of $L,$ hence $x_0=\frac{\beta}{2}z_0$ and $ H(x)<0$ if and only if $L(z)<0$.
 For $\theta > p \geq 1$, there holds
\begin{align*}
L(p+1) & = (p+1)^4 -\frac{16p\theta(p+1)^3}{(\theta+1)} +\frac{16p\theta(p+1)^2(p+\theta+2)}{(\theta+1)^2} -\frac{16p\theta(p+1)^2}{(\theta+1)^2}\\
& = (p+1)^4 -\frac{16p\theta(p+1)^3}{(\theta+1)} + \frac{16p\theta(p+1)^2}{(\theta+1)} + \frac{16p\theta(p+1)^3}{(\theta+1)^2}
-\frac{16p\theta(p+1)^2}{(\theta+1)^2}\\
& = (p+1)^2\left[(p+1)^2 -\frac{16p^{2}\theta}{(\theta+1)} + \frac{16p^{2}\theta}{(\theta+1)^2}\right]\\
&=\left(\frac{p+1}{\theta+1}\right)^{2}\left[(p+1)^2(\theta+1)^2 -16p^{2}\theta^{2}\right] < 0.
\end{align*}
The last inequality holds true since
$$4p\theta- (p+1)(\theta + 1) > 4p^2 - (p+1)^2 \geq 0, \quad \forall\; \theta > p \geq 1.$$
As $\lim_{z\rightarrow \infty}L(z)= \infty,$  it follows that $z_0>p+1.$ We get then
 \begin{align*}
2x_0>(p+1)\beta=2+\alpha+\beta,\quad \forall\; \theta > p\geq 1.
\end{align*}
If $p>\frac{4}{3}$, by $(i)$ of Theorem \textbf{B}, the system \eqref{1.1} has no classical stable solution if  $N_{s}<2+\alpha+\beta$. Suppose now $1\leq p \leq \min(\frac{4}{3}, \theta)$.
Observe that for all $\theta \geq p \geq 1$,
     \begin{align*}
\left[p+\frac{2(2-p)(p \theta -1)}{(\theta+p-2)(\theta+1)}\right]\beta \geq \alpha + \beta & \Leftrightarrow \left[p+\frac{2(2-p)(p \theta -1)}{(\theta+p-2)(\theta+1)}\right](\theta + 1) \geq p + \theta + 2\\
& \Leftrightarrow p\theta - 1 + \frac{2(2-p)(p \theta -1)}{\theta+p-2} \geq \theta + 1\\
& \Leftrightarrow (p\theta - 1)\left[ 1 + \frac{2(2-p)}{\theta+p-2}\right] \geq \theta + 1\\
& \Leftrightarrow (p\theta - 1)(\theta + 2 - p) \geq (\theta+p-2)(\theta + 1)\\
& \Leftrightarrow p\theta^2 - \theta + (2 - p)p\theta \geq \theta^2 + (p-1)\theta\\
& \Leftrightarrow (p-1)(\theta - p) \geq 0.
\end{align*}
As $z_0>p+1\geq 2,$ we have $x_0 = \frac{\beta z_0}{2} \geq \beta$ and
$$2+\alpha+\beta \leq 2 + \beta \left[p+\frac{2(2-p)(p \theta -1)}{(\theta+p-2)(\theta+1)}\right] \leq 2 + x_0\left[p+\frac{2(2-p)(p \theta -1)}{(\theta+p-2)(\theta+1)}\right].$$
If $N_{s} < 2+\alpha+\beta,$ using $(ii)$ of Theorem \textbf{B}, we are done.
\medskip
To conclude

\begin{align*}
\mbox{ for all} \quad \theta > p \geq 1, \mbox{ and} \quad N_{s} < 2+\alpha+\beta, \quad \eqref{1.1}\, \mbox{has no smooth stable solution}.
\end{align*}
\smallskip

The proof is finished.  \qed

\end{document}